\documentclass[12pt]{amsart}
\usepackage{amssymb}

\newtheorem{thm}{Theorem}[section]
\newtheorem{lem}[thm]{Lemma}
\newtheorem{cor}[thm]{Corollary}
\newtheorem{prop}[thm]{Proposition}

\newtheorem*{prob*}{Open problem}

\theoremstyle{definition}

\newtheorem{defi}[thm]{Definition}

\theoremstyle{remark}

\newtheorem{rem}[thm]{Remark}
\newtheorem*{rem*}{Remark}

\DeclareMathOperator{\id}{id}

\DeclareMathOperator{\Hom}{Hom}
\DeclareMathOperator{\Aff}{Aff}
\DeclareMathOperator{\Aut}{Aut}
\DeclareMathOperator{\Lie}{Lie}

\newcommand{\kringel}{\mathbin{\raise1pt\hbox{$\scriptstyle\circ$}}}
\newcommand{\pkt}{\mathbin{\raise0pt\hbox{$\scriptstyle\bullet$}}}

\newcommand{\R}{\mathbb{R}}

\newcommand{\tr}{\mathop{\rm tr}}
\newcommand{\ad}{\mathop{\rm ad}}

\newcommand{\End}{\mathop{\rm End}}
\newcommand{\Der}{\mathop{\rm Der}}

\newcommand{\La}{\mathfrak{a}}

\newcommand{\Lg}{\mathfrak{g}}

\newcommand{\Ll}{\mathfrak{l}}
\newcommand{\Ln}{\mathfrak{n}}

\newcommand{\CS}{\mathcal{S}}

\newcommand{\im}{\mathop{\rm im}}

\newcommand{\al}{\alpha}

\newcommand{\ov}{\overline}

\newcommand{\ra}{\rightarrow}

\renewcommand{\phi}{\varphi}

\begin{document}

\title[LR-algebras]{Complete LR-structures on solvable Lie algebras}

\author[D. Burde]{Dietrich Burde}
\author[K. Dekimpe]{Karel Dekimpe}
\author[K. Vercammen]{Kim Vercammen}
\address{Fakult\"at f\"ur Mathematik\\
Universit\"at Wien\\
  Nordbergstr. 15\\
  1090 Wien \\
  Austria}
\email{dietrich.burde@univie.ac.at}
\address{Katholieke Universiteit Leuven\\
Campus Kortrijk\\
8500 Kortrijk\\
Belgium}
\email{karel.dekimpe@kuleuven-kortrijk.be}
\email{kim.vercammen@kuleuven-kortrijk.be}

\date{\today}

\thanks{The first author thanks the KU Leuven Campus Kortrijk
for its hospitality and support}
\thanks{The second author expresses his gratitude towards the Erwin Schr\"odinger
International Institute for Mathematical Physics}
\thanks{K.\ Vercammen is supported by a Ph.~D.~fellowship of the Research Foundation - Flanders (FWO)}

\begin{abstract}
An LR-structure on a Lie algebra $\Lg$ is a bilinear product on $\Lg$, satisfying
certain commutativity relations, and which is compatible with the Lie product.
LR-structures arise in the study of simply transitive affine actions on Lie groups.
In particular one is interested in the question which Lie algebras admit a {\em complete}
LR-structure. In this paper we show that a Lie algebra admits a complete LR-structure
if and only if it admits any LR-structure.
\end{abstract}

\maketitle

\section{Introduction}

A Lie admissible algebra $(A,\cdot)$ is an algebra with a bilinear product $x\cdot y$
such that the commutator $[x,y]=x\cdot y-y\cdot y$ defines a Lie bracket.
Several Lie admissible algebra structures arise in geometry in a very natural way.
Particular examples are so called {\it left symmetric structures} and {\it LR-structures},
appearing in the study of left-invariant affine structures on Lie groups,
and simply transitive affine actions of Lie groups. There exists a large literature,
we refer to \cite{BEN, BUR, DEK, HEL, KIM, MED, MIL, MIZ, SEG} and the references given
therein.

\begin{defi}
A vector space $A$ over a field $k$ together with a bilinear product
$A\times A \rightarrow A$, $(x,y)\mapsto x\cdot y$ is called a
{\it left symmetric algebra} $(A,\cdot )$ if
\begin{align}
x\cdot(y\cdot z) - (x \cdot y) \cdot z = y\cdot(x\cdot z) - (y \cdot x) \cdot z
\end{align}
for all $x,y,z\in A$. An algebra $(A,\cdot)$ is called an {\it LR-algebra}, if the product
satisfies the identities
\begin{align}
x\cdot (y\cdot z)& = y\cdot (x\cdot z) \label{lr1}\\
(x\cdot y)\cdot z& =(x\cdot z)\cdot y \label{lr2}
\end{align}
for all $x,y,z \in A$.
\end{defi}

It is well known that every left symmetric algebra, and every LR-algebra is a Lie admissible
algebra, see \cite{BDD}, with associated Lie algebra denoted by $\Lg_A$.
Conversely, given a Lie algebra $\Lg$, it is an important question, whether there exists
a left-symmetric algebra (resp. an LR-algebra)
$(A,\cdot)$, such that $x\cdot y-y\cdot y$ coincides with the Lie bracket $[x,y]$ of $\Lg$,
i.e., with $\Lg_A=\Lg$. If so, it is said that $\Lg$ {\it admits} a left symmetric structure,
resp. an LR-structure. This algebraic existence question usually is a very deep question
in connection with the existence of certain geometric structures on Lie groups.
For example, the existence of left-invariant affine structures on Lie groups is directly
related to the existence of left symmetric structures on the corresponding Lie algebras,
and this is a very hard question for solvable Lie algebras, see \cite{MIL}.
For a long time, it was thought that all solvable real Lie algebras
would admit a left symmetric structure. However, in $1995$ Y.~Benoist \cite{BEN} gave
the first example of a nilpotent Lie algebra in dimension $11$ not admitting any
left symmetric structure. Later easier examples, and families of examples were constructed,
see \cite{BG, BUR}. Also, the completeness of left invariant affine structures on Lie
groups was studied. As left-invariant affine structures correspond to torsion-free, flat
connections $\nabla$ in the tangent bundle, completeness means that the corresponding
connection $\nabla$ is complete.
It is known that a left-invariant affine structure on a Lie group $G$ is complete if and
only if the right multiplications $R(x)$ in the corresponding left symmetric algebra $(A,\cdot)$ of
$\Lg=\Lie (G)$ are nilpotent for all $x\in \Lg$, see \cite{SEG}. Here the linear operators $R(x)$
are defined by $R(x)y=y\cdot x$. The following definition therefore makes sense.

\begin{defi}
A left symmetric algebra structure $(A,\cdot)$ on a Lie algebra $\Lg$ is said to be
{\it complete} if the right multiplication operators $R(x)$ are nilpotent for
all $x\in A$.
\end{defi}

It is known that a given connected and simply connected Lie group $G$ admits
a representation $\rho \colon G \rightarrow \Aff(\R^n)$ letting $G$ act simply transitively
on $\R^n$ if and only if the Lie algebra $\Lg$ of $G$ admits a complete left symmetric structure, see \cite{FG, FGH, KIM2}. In particular, the question of which Lie groups $G$ admit a
simply transitive and affine action  reduces to the question of which Lie algebras $\Lg$
(over $\R$) admit a complete left symmetric structure.
By a result of L. Auslander, such a Lie group $G$, and hence also the
Lie algebra $\Lg$, has to be solvable.
There was the hope, that conversely every connected and simply connected solvable
Lie group would admit a complete left symmetric structure. As mentioned above, this
was answered later in the negative. But before the counter examples it was also asked whether
a complete left symmetric structure would exist automatically, once there was any such
structure at all. Mizuhara showed that in the nilpotent case this is indeed the case:

\begin{thm}\cite{MIZ}\label{resultMIZ}
Let $\Lg$ be a complex nilpotent Lie algebra. If $\Lg$ admits a left symmetric structure,
then $\Lg$ also admits a complete left symmetric structure.
\end{thm}

Because of the counter examples mentioned, one has broadened the geometric context and
studied  affine actions on simply connected nilpotent Lie groups: let $N$ be such a connected
and simply connected nilpotent Lie group, then the affine group of $N$ is defined as the
semi-direct product $\Aff(N)=N \rtimes \Aut(N)$, where $\Aff(N)$ acts on $N$ via
\[ \forall m,n \in N,\; \forall \alpha \in \Aut(N):\; (m,\alpha)\cdot n= m \alpha(n).\]
It is easy to see that in case $N=\R^n$, the unique $n$-dimensional connected and simply connected
abelian Lie group, this concept of affine group coincides with the usual one. \\
In contrast to the fact that there are simply connected connected solvable Lie groups
$G$ not allowing a simply transitive affine action, all of them do admit a simply transitive
affine action on a nilpotent Lie group.

\begin{thm}(\cite{BAU},\cite{DEK})
Let $G$ be a connected and simply connected solvable Lie group, then there exists a
connected and simply connected nilpotent Lie group $N$ and a representation
$\rho \colon G\rightarrow \Aff(N)$ letting $G$ act simply transitively on $N$.
\end{thm}

In order to get a better understanding of such simply transitive affine actions $\rho:G \rightarrow
\Aff(N)$, we started in \cite{BDD1} and \cite{BDD} to investigate the situation
for $G=\R^n$, an abelian connected and simply connected Lie group, i.e.\ we have been
considering the case of abelian simply transitive subgroups of $\Aff(N)$ for a given
nilpotent Lie group. Also in this case a translation  to the Lie algebra level is
possible, and leads to LR-structures defined above.
Just as for left symmetric structures, we define the notion of complete LR-structures.

\begin{defi}\label{complete}
An LR-structure $(A,\cdot)$ on a Lie algebra $\Lg$ is called {\it complete}, if all
right multiplications $R(x)$ are nilpotent.
\end{defi}

\begin{rem}\label{opp}
We note that $\Lg$ admits an LR-structure such that all right multiplications
are nilpotent, if and only if $\Lg$ admits an LR-structure for which all left
multiplications are nilpotent. Indeed, let $x\cdot y$ denote a complete LR-structure,
then the opposite product $x\circ y:=-y\cdot x$ defines an LR-structure where all
left multiplications are nilpotent.
\end{rem}

As already announced above, also these structures have a geometrical meaning:

\begin{thm} Let $N$ be a connected and simply connected nilpotent Lie group with
Lie algebra $\Ln$. Then $N$ admits a simply transtive abelian action
$\rho \colon \R^n\rightarrow \Aff(N)$ if and only $\Ln$ admits a complete LR-structure.
\end{thm}

In \cite{BDD} we already obtained some general properties of LR-algebras. One of these results says
that a Lie algebra $\Lg$ admitting an LR-structure has to be 2-step solvable.
However, not every 2-step solvable Lie algebra admits an LR-structure.

\medskip

One goal of this paper is to obtain the analogue of Mizuhara's result for LR-algebras.
We will show that if a Lie algebra $\Lg$ admits any LR-structure, then $\Lg$ also
admits a complete LR-structure.

\medskip

Although the definitions introduced above also make sense in the infinite-dimensional
case, we will always assume in this paper, even if we do not mention this explicitly,
that all algebras are finite-dimensional, and defined over a field $k$ of characteristic
zero.

\medskip

{\bf Acknowledgement:} The authors thank Nansen Petrosyan for his helpful suggestions.

\section{LR-algebras}

In an LR-algebra several commutation rules hold.
If we denote by $L(x), R(x)$ the left respectively right multiplication operator
then $\eqref{lr1}$ and $\eqref{lr2}$ are equivalent to the requirement
that all left and all right multiplications commute:
\begin{align*}
[L(x),L(y)] & = [R(x),R(y)] =0.
\end{align*}

We can derive more commutation rules:

\begin{lem} \label{1.4}
Let $(A,\cdot)$ be an LR-algebra. Then the following identities hold in $A$:
\begin{align}
L(x)R(y) & =R(x\cdot y), \label{lr3a} \\
R(x)L(y) & =L(y\cdot x), \label{lr3b} \\
L(x)R(y\cdot z)& =R(x\cdot (y\cdot z)), \label{lr4} \\
R(x)L(y\cdot z)& =L((y\cdot z)\cdot x), \label{lr5} \\
L(x)L(y\cdot z)& =L(y\cdot (x\cdot z)), \label{lr6} \\
R(x)R(y\cdot z)& =R((y\cdot x )\cdot z) \label{lr7}.
\end{align}
\end{lem}

\begin{proof}
Let $x,y,z,a\in A$. Using \eqref{lr1} we have
\begin{align*}
L(x)R(y)(a)& = x\cdot (a\cdot y)\\
           & = a\cdot (x\cdot y)=R(x\cdot y)(a),
\end{align*}
so that $L(x)R(y)=R(x\cdot y)$. Rewriting $y$ as $y\cdot z$ we obtain \eqref{lr4}.
In the same way, using \eqref{lr2} we have
\begin{align*}
R(x)L(y)(a)& = (y\cdot a)\cdot x\\
           & = (y\cdot x)\cdot a=L(y\cdot x)(a),
\end{align*}
so that $R(x)L(y)=L(y\cdot x)$.  This implies \eqref{lr5}.
Furthermore, using \eqref{lr3a} and \eqref{lr3b}, we have
\begin{align*}
L(x)L(y\cdot z) & = L(x)R(z)L(y)\\
                & = R(x\cdot z)L(y)\\
                & = L(y\cdot (x\cdot z)),
\end{align*}
which shows \eqref{lr6}. In the same way follows \eqref{lr7}.
\end{proof}

\begin{prop}\label{nilp}
Let $(A,\cdot)$ be an LR-structure on a Lie algebra $\Lg$. Then any two of
the statements below imply the third one.
\begin{itemize}
\item[(a)] All left multiplications $L(x)$ are nilpotent operators.
\item[(b)] All right multiplications $R(x)$ are nilpotent operators.
\item[(c)] The Lie algebra $\Lg$ is nilpotent.
\end{itemize}
\end{prop}

\begin{proof}
Assume that $(a)$ and $(b)$ hold. We will show that $(c)$ holds.
For all $n\ge 1$ consider the identity
\[
\ad (x)^n = (L(x)-R(x))^n.
\]
By assumption all left and right multiplications, and their powers, have zero trace.
We claim that the right hand side of the above identity has zero trace for all $n\ge 1$.
For this it is enough to show that all summands between $L(x)^n$ and $R(x)^n$
are of the form $L(y\cdot z)$ or $R(y\cdot z)$, where $y,z$ are certain powers
of $x$, with various bracketings. We prove this by induction on $n$. The case
$n=1$ is clear. Assume this for $n$. Then
$\ad (x)^{n+1}=(L(x)-R(x))(L(x)-R(x))^n$.
The inner summands are of the form
\[
L(x)R(x)^n,\; R(x)L(x)^n,\;L(x)L(y\cdot z),\; L(x)R(y\cdot z),\;R(x)R(y\cdot z),\;
R(x)L(y\cdot z).
\]
They are all of the form  $L(y\cdot z)$ or $R(y\cdot z)$.
This follows from lemma $\ref{1.4}$.
For the last four summands this is obvious, and for the
first two summands we obtain
\begin{align*}
L(x)R(x)^n & =R(y\cdot x),\\
R(x)L(x)^n & =L(x\cdot z)
\end{align*}
for certain $y,z\in A$. For example, $L(x)R(x)=R(x\cdot x)$ and
$R(x)L(x)=L(x\cdot x)$.
This can be seen again by induction on $n$ and lemma $\ref{1.4}$.
It follows that $\tr(\ad(x)^n)=0$ for all $n\ge 1$, so that all $\ad(x)$ are nilpotent,
because the field has characteristic zero. By Engel's theorem $\Lg$ is nilpotent. \\[0.2cm]
Assume that $(b)$ and $(c)$ hold. We will prove that $(a)$ holds.
Now we consider the identity
\[
L(x)^n=(\ad(x)+R(x))^n,
\]
where $n\ge 1$ and $x\in \Lg$. We claim that the right hand side of this identity has
zero trace for all $n\ge 1$. This follows because it can be expressed as a linear
combination of the terms
\[
\ad(x)^n,R(x)^n,L(y\cdot z), R(y\cdot z),
\]
where $y$ and $z$ are certain powers of $x$ with various bracketings.
To verify this, one does a very similar calculation as in the previous case.
By assumption all $\ad(x)$ and all $R(x)$ are nilpotent. Because
$L(x)=\ad(x)+R(x)$, all $L(x)$ have zero trace. Then all terms listed above have
zero trace, so that $L(x)^n$ has zero trace for all $n\ge 1$, and we are done. \\[0.2cm]
Finally, assume that $(a)$ and $(c)$ hold. We have to show that $(b)$ holds.
This case follows from the above case using remark $\ref{opp}$.
\end{proof}

\section{Nilpotent Lie algebras}

For the study of LR-structures on $\Lg$ we repeatedly need the following proposition, which can
be found in \cite[Proposition 5.3]{BD}:

\begin{prop}\label{module-split}
Let $\Lg$ be a nilpotent Lie algebra over a field $k$ of characteristic zero,
and let $V$ be a finite dimensional $\Lg$-module. Then $V$ can be written as a
direct sum of $\Lg$-modules $V=V_n\oplus V_0$ where $V_n$ is the unique maximal
nilpotent submodule of $V$ and $H^i(\Lg,V_0)=0,\; \forall i\geq 0$.
\end{prop}

Our first main result is the following:

\begin{prop}\label{mizuhara}
Let $\Lg$ be a finite-dimensional nilpotent Lie algebra over $k$.
If $\Lg$ admits an LR-structure then it also admits a complete LR-structure.
\end{prop}

\begin{proof}
Let $(A,\cdot)$ be an LR-structure on $\Lg$. We will first assume that $k$ is algebraically
closed. Denote by $\La=\im (L)$ the subalgebra of $\Lg\Ll(\Lg)$ generated by all $L(x)$.
Since these operators all commute this subalgebra is abelian, in particular
nilpotent. Denote by
\[
\rho\colon \La \hookrightarrow \Lg\Ll(\Lg)
\]
the representation of $\La$ given by inclusion. Since $k$ is algebraically closed we can
apply the weight subspace decomposition for $\rho$. We have
\[
A=\bigoplus_{i=1}^s A_{\al_i}(\La),
\]
where $\al_i\in \Hom (\La,k)$ are the different weights of $\rho$ and
$A_{\al_i}(\La)=A_i$ are the (nonzero) weight subspaces.
We may identify  $\rho(L(x))$ with $L(x)$. In a suitable basis of the vector
space $A$ the operators $L(x)$ have a block matrix, where each block is of the
form
\[
L(x)_{\mid A_i}=
\begin{pmatrix}
\al_i(L(x)) &  &  & \\
         &  &  & \ast \\
         & \al_i(L(x)) & & \\
0        &  & \ddots & \\
         &  &        &  \al_i(L(x))\\
\end{pmatrix}.
\]
A first consequence is that each $A_i$ is a left ideal in $A=A_1\oplus \cdots \oplus A_s$.
For each $1\le i \le s$ let $(e_{i,1},\ldots ,e_{i,n_i})$ be a basis of $A_i$.
Then $\{e_{i,1},\ldots ,e_{i,n_i}\}_{i=1,\ldots ,s}$ is a suitable basis for $A$
as mentioned above. Let $A_{i,k}=\langle e_{i,1},\ldots e_{i,k}\rangle$ for $k=1,\ldots ,n_i$.
We obtain a  a filtration of $A_i$ by left ideals:
\[
A_i=A_{i,n_i}\supseteq A_{i,n_i-1}\supseteq \cdots \supseteq A_{i,2} \supseteq A_{i,1}
\supseteq A_{i,0}=0.
\]
We can write each $x\in A$ uniquely as $x=x_1+\cdots +x_s$ with $x_i\in A_i$, i.e., as
\[
x=\sum_{i=1}^s (x_{i,1}e_{i,1}+\ldots + x_{i,n_i}e_{i,n_i}).
\]
\begin{lem}\label{2.2}
Suppose that we have in the representation of $x\in A$ that $x_i\in A_{i,n_i-1}$ for some $i$,
i.e., $x_{i,n_i}=0$. Then $\al_i(L(x))=0$.
\end{lem}

\begin{proof}
Assume that $\al_i(L(x))\neq 0$. We will show that this is impossible.
For $z\in A$ the condition $x_i\in A_{i,n_i-1}$ implies that $L(z)x_i \in A_{i,n_i-1}$, because
$A_{i,n_i-1}$ is a left ideal. This means that $(L(z)x)_{i,n_i}=0$ for all $z\in A$. \\
We can choose an element $y\in A$ satisfying
$y_{i,n_i}=1$ with respect to its basis representation. Consider the elements
$y^{(0)}=y$ and $y^{(k)}=[x,y^{(k-1)}]$ for all $k\ge 1$. Using induction it is
easy to see that $(y^{(k)})_{i,n_i}=(\al_i(L(x)))^k$. Indeed,
\begin{align*}
(y^{(k+1)})_{i,n_i} & = (L(x)y^{(k)})_{i,n_i}-(L(y^{(k)})x)_{i,n_i}\\
                   & = (L(x)y^{(k)})_{i,n_i} \\
                   & =  \al_i(L(x))\, (y^{(k)})_{i,n_i} \\
                   & = (\al_i(L(x)))^{k+1}.
\end{align*}
This shows that all elements $y^{(k)}$ are nonzero. In particular, the operator $\ad (x)$
for this element $x\in A$ is not nilpotent. Hence the Lie algebra $\Lg$ is not nilpotent,
which is a contradiction.
\end{proof}

{\it Case 1: One weight.} We first assume now that $\rho$ has only a single nonzero
weight $\al$. Then there is a basis $(e_1,\ldots ,e_n)$ of $A$ such that
the operators $L(x)$ are upper-triangular with $\al(L(x))$ on the diagonal.
For $x\in A$ we write $x=x_1e_1+\cdots +x_ne_n$. Lemma $\ref{2.2}$ says that $\al(L(x))=0$
for elements $x$ with $x_n=0$, i.e., $\al(L(e_k))=0$ for $k=1,\ldots ,n-1$. Since $\al$
is nonzero it follows that $\al(L(e_n))\neq 0$, and we may normalize it to $1$.
We may write the operators $L(e_k)$ as follows:
\[
L(e_k)=(a_{i,j}^k)_{1\le i,j \le n}, \quad k=1,\ldots ,n,
\]
with the conditions
\begin{align*}
a_{i,j}^k & = 0 \text{ for all } 1\le j\le i \le n, 1\le k\le n-1, \\
a_{i,j}^n & = 0 \text{ for all } 1\le j<i \le n, \\
a_{i,i}^n & = 1 \text{ for all } 1\le i\le n.
\end{align*}

We want to show that $A$ in this case is commutative. For this we need two lemmas.

\begin{lem}\label{2.3}
We have $a_{n-1,n}^k=0$ for all $k=1,\ldots ,n-2$. In particular, the subspace
generated by $e_1,\ldots ,e_{n-2}$ is a two-sided ideal in $A$.
\end{lem}
\begin{proof}
The operators $R(e_j)$ for $1\le j\le n$ are given as follows:
\[
R(e_j)=
\begin{pmatrix}
a_{1,j}^1 & a_{1,j}^2 & \cdots & a_{1,j}^{n-1} & a_{1,j}^{n} \\
a_{2,j}^1 & a_{2,j}^2 & \cdots & a_{2,j}^{n-1} & a_{2,j}^{n} \\
\hdots   &  \hdots   & \hdots & \hdots       & \hdots     \\
a_{j-1,j}^1 & a_{j-1,j}^2 & \cdots & a_{j-1,j}^{n-1} & a_{j-1,j}^{n} \\
0          & 0          & \cdots & 0              & 1 \\
0          & 0          & \cdots & 0              & 0 \\
\hdots   &  \hdots   & \hdots & \hdots       & \hdots     \\
0          & 0          & \cdots & 0              & 0 \\
\end{pmatrix}
\]
Suppose that it is not true that $a_{n-1,n}^1=a_{n-1,n}^2=\cdots = a_{n-1,n}^{n-2}=0$.
Then there is a minimal $k\le n-2$ such that $a_{n-1,n}^k\neq 0$, and
$a_{n-1,n}^j=0$ for all $1\le j\le k-1$. Then we consider the matrix identity
\[
[R(e_k),R(e_n)]=0.
\]
Looking at the $(n-1,n)$-th entry of the left hand side, we find that
 $a_{n-1,n}^k=0$, which is a contradiction. Hence we have
$a_{n-1,n}^k=0$ for all $k=1,\ldots ,n-2$. This implies that
$R(e_n)$ maps $\langle e_1,\ldots ,e_{n-2}\rangle$
into itself. This, together with the particular form of the operators
$R(e_j)$ now shows that $\langle e_1,\ldots ,e_{n-2}\rangle$ is also a right ideal
in $A$, which finishes the proof.
\end{proof}

\begin{lem}\label{2.4}
There is an idempotent $a\neq 0$ in $A$, i.e., satisfying $a\cdot a=a$.
\end{lem}

\begin{proof}
We will use induction on $n$. For $n=1$ we have
$L(e_1)=(1)$, so that $e_1\cdot e_1=e_1$. For $n=2$ the left multiplications have the
form
\[
L(e_1)=\begin{pmatrix}
0 & a_{1,2}^1 \\
0 & 0
\end{pmatrix}, \quad
L(e_2)=\begin{pmatrix}
1 & a_{1,2}^2 \\
0 & 1
\end{pmatrix}.
\]
We have $a_{1,2}^1=1$. Indeed, $\tr(\ad (e_2))=\tr(L(e_2))-\tr(R(e_2))=1-a_{1,2}^1$ has to be
zero. Now $a:=-a_{1,2}^2e_1+e_2$ is the desired non-trivial idempotent. \\
Suppose $n\ge 3$. By lemma $\ref{2.3}$ the subspace
$B=\langle e_1,\ldots ,e_{n-2}\rangle$ is an ideal in $A$ such that $A/B$ is an
LR-structure on the associated Lie algebra $\Lg_{A/B}$, which is again nilpotent.
The associated homomorphism $\ov{\rho}$ again has one single nonzero weight. Since
$\dim (A/B)=2$ there is a non-trivial idempotent $\ov{a}=a+B$ in $A/B$. We have
$a\cdot a+B=\ov{a}\cdot \ov{a}=\ov{a}$, hence $a\cdot a=a+b$ for some $b\in B$.
Note that $a\not\in B$. For such an element $a$ let $B_a=\langle B,a\rangle$.
This is an LR-algebra with associated nilpotent Lie algebra and
homomorphism having again one single nonzero weight. Since
$\dim B_a<\dim A$ we can apply the induction hypothesis: there is a non-trivial
idempotent $e\in B_a$, hence also in $A$.
\end{proof}

Now we can show that $A$ is commutative:

\begin{prop}\label{2.5}
If $\rho$ has only a single nonzero weight $\al$, then $A$ is commutative.
Hence the zero product defines a complete LR-structure on $\Lg$.
\end{prop}

\begin{proof}
Let $a\neq 0$ be an idempotent in $A$. Then $L(a)$ is not nilpotent. Because
$L(a)$ is upper-triangular with nonzero diagonal elements $\al(L(a))$,
the operator $L(a)$ is invertible. Lemma $\ref{1.4}$ implies that
$R(a)L(a)=L(a\cdot a)=L(a)$.
Multiplying with $L(a)^{-1}$ yields $R(a)=\id$. This means, $a$ is a right identity
for $A$. Again using lemma  $\ref{1.4}$ we obtain
\begin{align*}
L(x)=L(x)R(a)=R(x\cdot a)=R(x),
\end{align*}
so that $A$ is commutative, and $\Lg$ is abelian.
\end{proof}
\vspace*{0.5cm}
{\it Case 2:} We now consider the general case, i.e., where $A=A_1\oplus \cdots \oplus A_s$
is the direct sum of weight spaces $A_i$ corresponding to the weights $\al_i$.
We may assume that all $\al_i$ with $i\ge 2$ are nonzero, and that $\al_1=0$, possibly
setting $A_1=0$ if there is no zero weight.
The spaces $A_i$ for $i\ge 2$ are left ideals, hence LR-algebras with associated
nilpotent Lie algebra and homomorphism $\rho_{\mid L(A_i)}$ having one single weight.
It is a consequence of lemma $\ref{2.2}$ that this weight is nonzero.
Now it follows from the first case that all $A_i$ for $i\ge 2$ are commutative
and have non-trivial right identities $a_i\in A_i$. In particular, if $x\in A_i$
and $y\in A_j$ for $i\neq j$ with $2\le i\le s$ and $1\le j\le s$ we have
\begin{align}\label{8}
x\cdot y & =(x\cdot a_i)\cdot y =(x\cdot y)\cdot a_i \; \in A_i\cap A_j\; =0.
\end{align}

Now define a bilinear product $x\circ y$ on $A$ as follows. On basis vectors we set
\[
e_{i,k}\circ e_{j,l}= \begin{cases} 0 \quad \hspace{1cm} \forall\; i=j\ge 2,\;
1\le k,l\le n_i  \\
e_{i,k}\cdot e_{j,l} \; \text{ in all other cases.} \end{cases}
\]
The proof is finished if we can show that this product defines a
complete LR-structure on $\Lg$.
Note that for $1\le k\le n_i$ and $1\le l\le n_j$ we have
\begin{align}
e_{i,k}\circ e_{j,l} & = 0 \quad \forall \, i\ge 2,  \; \forall j\ge 1, \label{st1}\\
[e_{i,k},e_{j,l}]    & = 0 \quad \forall i,j\ge 2. \label{st2}
\end{align}
From this it follows easily that
\[
[e_{i,k},e_{j,l}] = e_{i,k}\circ e_{j,l} - e_{j,l}\circ e_{i,k}.
\]
Indeed, if $i,j\ge 2$ then both sides are equal to zero, while in the
other cases the product $x\circ y$ coincides with the original LR-product $x\cdot y$. \\
Next we will show that all left multiplications for this product commute, i.e., that
\[
e_{i,p}\circ (e_{j,q}\circ e_{k,r})=e_{j,q}\circ (e_{i,p}\circ e_{k,r}).
\]
If $j\ge 2$ or $i\ge 2$ then both sides are equal to zero by \eqref{st1}.
In the other cases the product $x\circ y$ coincides with the
original LR-product $x\cdot y$, which satisfies this identity. \\
Next we show that all right multiplications commute, i.e., that
\[
(e_{i,p}\circ e_{j,q})\circ e_{k,r}=(e_{i,p}\circ e_{k,r})\circ e_{j,q}.
\]
For $i\ge 2$ both sides are equal to zero by \eqref{st1}. Thus we can assume that
$i=1$. Then if both $j\ge 2$ and $k\ge 2$ we obtain $0=0$, because $A_j$ and $A_k$
are left ideals, and we can apply \eqref{st1}. In the case $j=k=1$ we can replace
the product $x\circ y$ by $x\cdot y$, and the identity is satisfied.
It remains to check the cases with $i=1$ and either $j$ or $k$ equal to $1$. By
symmetry we may assume that $j=1$ and $k\ge 2$. Then the identity to show
is
\[
(e_{1,p}\cdot e_{1,q})\circ e_{k,r}=(e_{1,p}\cdot e_{k,r})\circ e_{1,q}.
\]
The left hand side equals $(e_{1,p}\cdot e_{1,q})\cdot e_{k,r}$, because $A_1$
is a left ideal. The right hand side is zero, because $A_k$ is a left ideal.
But by \eqref{8} we have $0=(e_{1,p}\cdot e_{k,r})\cdot e_{1,q}$, so that
the above identity reduces to
\[
(e_{1,p}\cdot e_{1,q})\cdot e_{k,r}=(e_{1,p}\cdot e_{k,r})\cdot e_{1,q},
\]
which holds true.
Finally we prove that the new LR-structure is complete.
Denote by $\ell(x)$ the left multiplications, i.e., $\ell(x)y=x\circ y$,
and by $r(x)$ the right multiplications.
We see that all $\ell(e_{i,p})$ are nilpotent, because each of its restrictions
to one of the left ideals $A_j$ is nilpotent. Indeed, for $j=1$ we have
$\ell(e_{i,p})_{\mid A_1}=L(e_{i,p})_{\mid A_1}$, which is nilpotent since $A_1$ is the
weight space corresponding to the weight zero. If $j\ge 2$ then we have
$\ell(e_{i,p})_{\mid A_j}=0$ for $i\ge 2$. If $i=1$ then $\ell(e_{1,p})_{\mid A_j}=L(e_{1,p})_{\mid A_j}$
is nilpotent by lemma $\ref{2.2}$.
Because all $\ell(e_{i,p})$ commute, we obtain
that all $\ell(x),x\in A$ are nilpotent. Since $\Lg$ is nilpotent,
proposition $\ref{nilp}$ implies that also all right multiplications $r(x)$
are nilpotent. Hence the new LR-structure is complete. \\

Finally let us consider the case where $\Lg$ is a nilpotent Lie algebra over
an arbitrary field $k$ of charactersitic 0. It is still true that $\La=\im(L)$
is an abelian Lie algebra and that $A$ is an $\La$-module. By proposition~\ref{module-split}
we know that $A$ splits as a direct sum of $\La$-modules: $A=V_n\oplus V_0$ where $V_n$ is the
unique maximal nilpotent submodule of $A$. Stated differently, $V_n$ is the unique
maximal left ideal of $A$ on which $A$ acts nilpotently by multiplications from the left.

If we denote the algebraic closure of $k$, by $\bar{k}$, then $\Lg\otimes \bar{k}$ is a
Lie algebra over $\bar{k}$ having an LR-structure $(A\otimes \bar{k},\cdot)$.
Note that $V_n\otimes \bar{k}$ is the unique maximal left ideal of  $A\otimes
\bar{k}$ on which $A\otimes \bar{k}$ acts nilpotently by multiplications from the left.

It follows that when we decompose $A\otimes \bar{k}$ as a direct sum $A_1\oplus A_2 \oplus \cdots
\oplus A_s$ of its weight spaces as we did before, then the subspace
corresponding to the weight zero is $A_1=V_n\otimes \bar{k}$ and
$A_2\oplus \cdots \oplus A_s=V_0\otimes \bar{k}$. So
any element $x\in A\otimes \bar{k}$ can be uniquely written as a sum
$x= x_n + x_0$, with $x_n\in  V_n\otimes \bar{k}$ and $x_0 \in V_0\otimes \bar{k}$.
Moreover, such an element belongs to $A=A\otimes k\subseteq A\otimes \bar{k}$
if and only if $x_n\in V_n=V_n\otimes k$ and $x_0\in V_0=V_0\otimes k$.

Now, let $\circ$ again denote the complete LR-structure on $\Lg\otimes \bar{k}$ as
constructed above. Using the definition of $\circ$ and
\eqref{st1}, we see that
\[ \forall x_n,y_n \in  V_n\otimes \bar{k},\;
   \forall x_0,y_0 \in  V_0\otimes \bar{k}:\;
   (x_n + x_0) \circ (y_n + y_0) = x_n \cdot y_n + x_n \cdot y_0 .\]
From this, it is obvious that $\circ$ restricts to a $k$-bilinear product
on $A=A\otimes k$, which is then of course a complete LR-structure on $A$.
\end{proof}

\begin{rem}\label{relation-new-old}
For use later on, we note that the last part of the above proof shows that
the complete LR-product $x\circ y$ on $\Lg$, which was constructed from the
original LR-product $x\cdot y$ satisfies $\Lg \circ \Lg \subseteq \Lg\cdot \Lg$.
\end{rem}

\section{Solvable Lie algebras}

For any Lie algebra $\Lg$, we denote by $\Lg^1=\Lg$ and $\Lg^{i+1}=[\Lg,\Lg^i]$ for
$i\ge 1$ the ideals of the lower central series of $\Lg$. As we always assume that $\Lg$
is finite dimensional, this series stabilizes after finitely many steps, with
\[
\Lg^{\infty}=\bigcap_{i=1}^{\infty}\Lg^{i}.
\]
Let $\Ln:=\Lg/\Lg^{\infty}$. This is a nilpotent Lie algebra.

\begin{lem}
Let $\Lg$ be a two-step solvable Lie algebra. Then
the extension
\[
0\ra \Lg^{\infty} \ra \Lg \ra \Lg/\Lg^{\infty} \ra 0
\]
splits, so that $\Lg=\Lg^{\infty} \rtimes \Ln$.
\end{lem}

\begin{proof}
As $\Lg^{\infty}\subset [\Lg,\Lg]$,
it follows that $[\Lg^{\infty},\Lg^{\infty}]\subseteq  [[\Lg,\Lg],[\Lg,\Lg]]=0$, because
$\Lg$ is two-step solvable. So $\Lg^{\infty}$ is abelian. The above short exact sequence
induces a Lie algebra homomorphism
\[
\phi\colon \Ln\ra \Der (\Lg^{\infty})= \End (\Lg^{\infty}).
\]
Recall that for all $n\in \Ln$ and all $x\in \Lg^\infty$: $\phi(n)(x)=[\tilde{n},x]$, where
$\tilde{n}$ is any pre-image of $n$ in $\Lg$.

Since $\Ln$ is nilpotent, we can apply proposition~\ref{module-split}, where we view
$\Lg^{\infty}$ as a $\Ln$-module via the representation
$\phi$. Since $[\Lg,\Lg^{\infty}]=\Lg^{\infty}$, the unique maximal nilpotent
submodule of $\Lg^{\infty}$ is trivial and hence
\[
H^i(\Ln,\Lg^{\infty})=0.
\]
for all $i\ge 0$. For $i=2$ this
means that the above extension splits.
\end{proof}

We are now ready to prove the main theorem of this paper.

\begin{thm}\label{mainthm}
Let $\Lg$ be a Lie algebra over a field $k$ of characteristic 0. If $\Lg$ admits an
LR-structure, then $\Lg$ also admits a complete LR-structure.
\end{thm}

\begin{proof}
As $\Lg$ admits an LR-structure $x\cdot y$, we know that $\Lg$ is 2-step solvable, see
\cite[Proposition 2.1]{BDD}. By the lemma above, we have that $\Lg=\Lg^{\infty}\rtimes \Ln$,
where $\Ln=\Lg/\Lg^{\infty}$ is a nilpotent Lie algebra. Since all terms of the lower
central series of $\Lg$ are 2-sided ideals for the LR-product, see \cite[Corollary 2.9]{BDD},
$\Lg^{\infty}$ is a 2-sided ideal of $\Lg$ and hence the LR-structure on $\Lg$
induces an LR-structure on $\Ln=\Lg/\Lg^\infty$. Let us denote this induced product on $\Ln$ by
$x\bullet y$.
By proposition~\ref{mizuhara} and remark~\ref{relation-new-old}
we know that $\Ln$ also admits a complete LR-structure $x\circ y$ with
$\Ln \circ \Ln\subseteq \Ln\bullet \Ln$.
When we now view $\Ln$ as a subspace of $\Lg=\Lg^{\infty}\rtimes \Ln$, we get that
$\Ln\circ \Ln \subseteq \Ln\cdot \Ln +\Lg^\infty$.

As above, let $\phi:\Ln \rightarrow \Der(\Lg^\infty)$ denote the homomorphism induced
by the split-extension, i.e.\ $\phi(n)(x)= [n,x]$ for all $n\in \Ln$ and $x\in \Lg^\infty$.

We claim that $\phi(\Ln\circ \Ln)=0$. Indeed, since ,$\Lg^\infty$ is abelian,
and $\Lg^\infty\subseteq [\Lg,\Lg]\subseteq \Lg\cdot \Lg$, we have
\begin{eqnarray*}
[\Ln\circ \Ln,\Lg^\infty ] & \subseteq & [\Ln \cdot \Ln + \Lg^\infty, \Lg^\infty]  \\
& \subseteq & [\Ln \cdot \Ln,  \Lg^\infty]  \\
& \subseteq & [\Lg \cdot \Lg, \Lg \cdot \Lg] \\
& = & 0
\end{eqnarray*}
where the last equality follows from the fact that in any LR-algebra $A$, we have the identity
$(x\cdot y)\cdot (u\cdot v) - (u\cdot v)\cdot (x\cdot y)=0$ for all $x,y,u,v\in A$, see
the proof of Proposition~2.1 in \cite{BDD}. \\
It now follows from Corollary~5.2 of \cite{BDD}, that $\Lg$ admits a
complete LR-structure $\star$, which is given by the following formula:
\[\forall (a,x), (b,y) \in \Lg=\Lg^\infty \rtimes \Ln:\;(a,x)\star (b,y)
= (\phi(x) b, x \circ y).\]
\end{proof}

We already saw that if a solvable Lie algebra $\Lg$ admits an LR-structure,
then also $\Ln=\Lg/\Lg^{\infty}$ admits an LR-structure.
The proof above also suggests a partial converse of this.

\begin{thm}\label{3.2}
Let $\Lg$ be a two-step solvable Lie algebra and assume that $\Ln=\Lg/\Lg^{\infty}$
admits an LR-structure satisfying
\begin{align}\label{9}
\Ln \cdot \Ln \subseteq [\Ln,\Ln].
\end{align}
Then, this LR-structure can be lifted to an LR-structure on $\Lg$.
\end{thm}

\begin{proof} As before, $\Lg=\Lg^{\infty} \rtimes \Ln$, where the action of $\Ln$ on $\Lg^{\infty}$
is given by $\phi\colon \Ln\ra \End (\Lg^{\infty})$. By assumption
$[[\Ln,\Ln],\Lg^{\infty}]\subseteq [[\Lg,\Lg],[\Lg,\Lg]]=0$, so that $\phi (\Ln\cdot\Ln)=0$.
If $(\Ln,\cdot)$ denotes the LR-structure on $\Ln$, then the same formula
$(a,x)\star (b,y) = (\phi(x) b, x \cdot y )$ as before defines an LR-structure on $\Lg$.
\end{proof}

\begin{rem} By an induction argument, as in the proof of \cite[Lemma 2.10]{BDD}, one can show
that \eqref{9} implies
\[
\Ln^i\cdot \Ln^j\subseteq \Ln^{i+j}
\]
for all $i,j\ge 1$. From this, it follows easily that the LR-structure $\star$ defined
above is a complete structure.
\end{rem}

\begin{cor}
Let $\Lg$ be a two-step solvable Lie algebra with $\Lg^{\infty}=\Lg^3$. Then
$\Lg$ admits an LR-structure.
\end{cor}

\begin{proof}
By assumption $\Ln=\Lg/\Lg^{\infty}=\Lg/\Lg^3$ is two-step nilpotent.
By proposition $4.3$ of \cite{BDD}, $x\cdot y=\frac{1}{2}[x,y]$ defines an LR-structure
on $\Ln$. It satisfies $\Ln \cdot \Ln \subseteq [\Ln,\Ln]$, so the claim
follows from the above theorem.
\end{proof}

\begin{rem}
There are two-step solvable Lie algebras with $\Lg^{\infty}=\Lg^4$ without any
LR-structure, see proposition $4.7$ of \cite{BDD}. In this sense the corollary
cannot be improved.
\end{rem}

\section{Two-step solvable two-generated Lie algebras}

Among the two-step solvable Lie algebras with two generators are the
filiform nilpotent Lie algebras of solvability class $2$. For this class
of Lie algebras we have constructed an explicit LR-structure in
\cite{BDD}. We want to generalize this construction here to all two-generated,
two-step solvable Lie algebras. We need two lemmas.

\begin{lem}\label{4.1}
Let $\Lg$ be a two-step solvable Lie algebra, $x_i,y\in \Lg$ for
$1\le i\le n$ and $\sigma$ a permutation in $\CS_n$. Then
\[
\ad (x_1)\cdots \ad (x_n)\ad(y)=\ad(x_{\sigma(1)})\cdots \ad(x_{\sigma(n)})\ad (y).
\]
\end{lem}

\begin{proof}
We use induction on $n\ge 1$, where the case $n=1$ is clear. For $n\ge 2$ we consider
$\sigma\in \CS_n$ with $\sigma(1)=1$. Then by the induction hypothesis we have
\begin{align}\label{10}
\ad(x_1) \ad(x_{\sigma(2)})\cdots \ad(x_{\sigma(n)})\ad (y)
& =\ad (x_1)\, \ad(x_2) \cdots \ad (x_n)\ad(y).
\end{align}
For $z\in \Lg$ the Jacobi identity and the solvability class $2$ imply
\begin{align*}
\ad(x_1)\ad(x_2)\ad (y)(z) & =  [x_1,[x_2,[y,z]]]\\
                           &  =-[x_2,[[y,z],x_1]]-[[y,z],[x_1,x_2]] \\
                           & = [x_2,[x_1,[y,z]]] \\
                           & = \ad(x_2)\ad(x_1)\ad(y)(z).
\end{align*}
In the same way we obtain
\begin{align}\label{11}
\ad(x_1)\ad(x_2)\cdots \ad(x_n)\ad (y)(z) & = \ad(x_2)\ad(x_1) \cdots \ad (x_n)\ad(y)(z).
\end{align}
The general case follows from \eqref{10} and \eqref{11}.
\end{proof}

\begin{lem}
Let $\Lg$ be a two-step solvable Lie algebra and $x,y\in \Lg$. Then the subspace
spanned by $x$ and
\[
\{\ad (y)^k \ad(x)^{\ell}y\mid k,l\ge 0 \}
\]
is a Lie subalgebra of $\Lg$.
\end{lem}

\begin{proof}
The bracket of each two generators is again in this subspace. If $\ell\ge 1$ then
lemma $\ref{4.1}$ implies
\begin{align*}
[x,\ad (y)^k \ad(x)^{\ell}y] & = \ad(x) \ad (y)^k \ad(x)^{\ell}y \\
 & = \ad (y)^k \ad(x)^{\ell+1}y.
\end{align*}
For $\ell,n\ge 1$ we have
\[
[\ad (y)^k \ad(x)^{\ell}y, \ad (y)^m \ad(x)^{n}y]=0,
\]
because $\Lg$ is two-step solvable. The other cases are trivial.
\end{proof}

Now we can show the following result.

\begin{thm}\label{4.3}
Let $\Lg$ be a two-generated, two-step solvable Lie algebra. Then $\Lg$ admits a
complete LR-structure.
\end{thm}

\begin{proof}
Let $\Lg$ be generated as a Lie algebra by $x$ and $y$. Then $\Lg$ is spanned
by $x$ and all vectors $\ad (y)^k \ad(x)^{\ell}y$ with $k,\ell\ge 0$. Let us fix
a basis consisting of $x,y$ and a subset of the above vectors satisfying $\ell \ge 1$.
Using this basis we define a $k$-bilinear product on $\Lg$ as follows
\begin{align*}
L(x) & = 0,\\
L(\ad(y)^k\ad(x)^{\ell}y) & = \ad(y)^k\ad(x)^{\ell}\ad(y).
\end{align*}
In particular, $L(y)=\ad (y)$. We have to verify that this
defines an LR-structure on $\Lg$. First of all we have $a\cdot b-b\cdot a=[a,b]$
for all basis vectors $a,b\in \Lg$: for $a=x$ and $b=\ad (y)^k \ad(x)^{\ell}y$ this means
\begin{align*}
a\cdot b-b\cdot a & = L(a)b-L(b)a \\
 & = 0-\ad(y)^k\ad(x)^{\ell}\ad (y)x \\
 & = \ad(y)^k\ad(x)^{\ell+1}y \\
 & = \ad(x)b \\
 & = [a,b]
\end{align*}
by lemma $\ref{4.1}$. For $a=\ad(y)^k\ad(x)^{\ell}y$ and $b=\ad(y)^m\ad(x)^{n}y$
with $\ell,n\ge 1$ this means
\begin{align*}
a\cdot b-b\cdot a & = \ad(y)^k \ad(x)^{\ell}\ad(y)\ad(y)^m\ad(x)^ny \\
                  & \; -  \ad(y)^m \ad(x)^{n}\ad(y)\ad(y)^k\ad(x)^{\ell}y \\
 & = \ad(y)^{k+m+1}\ad(x)^{\ell+n}y-  \ad(y)^{m+k+1}\ad(x)^{\ell+n}y\\
 & = 0 \\
 & = [a,b],
\end{align*}
since $\Lg$ is two-step solvable. The other cases are similar. Next we verify
that $[L(a),L(b)]=0$ for all basis vectors $a,b\in \Lg$. For
$a=\ad(y)^k\ad(x)^{\ell}y$ and $b=\ad(y)^m\ad(x)^{n}y$ with $\ell,n\ge 1$ this means
\begin{align*}
L(a) L(b) & = \ad(y)^{k+m+1}\ad(x)^{\ell+n}\ad(y) \\
               & = L(b)L(a).
\end{align*}
The other cases are similar. Finally we have to show that $[R(a),R(b)]=0$
for all basis vectors $a,b\in \Lg$. We have $R(a)=L(a)-\ad (a)$, so that
\begin{align*}
R(x) & = -\ad (x),\\
R(y) & = 0.
\end{align*}
It is not difficult to see that, for $a=\ad(y)^k\ad(x)^{\ell}y$ with $\ell \ge 1$
we have
\[
R(a)=\ad(y)^{k+1}\ad(x)^{\ell}.
\]
For $a=\ad(y)^k\ad(x)^{\ell}y$ and $b=\ad(y)^m\ad(x)^{n}y$ with $\ell,n\ge 1$
we have
\begin{align*}
R(a) R(b) & = \ad(y)^{k+1}\ad(x)^{\ell} \ad(y)^{m+1}\ad(x)^{n}\\
               & = \ad (y)^{m+1}\ad(x)^n\ad(y)^{k+1}\ad (x)^{\ell} \\
               & = R(b) R(a).
\end{align*}
The other cases are similar. It follows that $\Lg$ admits an LR-structure.
By theorem $\ref{mainthm}$ it admits also a complete LR-structure.
\end{proof}

\end{document}